\input amstex
\loadbold
\documentstyle{amsppt}
\TagsOnRight

\topmatter

\title
On the structure of finite-sheeted coverings of compact connected groups
\endtitle

\rightheadtext{coverings of compact connected groups}

\author
S. A. Grigorian and R. N. Gumerov
\endauthor

\address
Institute of Mechanics and Mathematics, Kazan State University,
Kremlevskaya 17, \ Kazan, \ 420008, \ Tatarstan, \ Russian Federation
\endaddress

\email
suren.grigorian\@ksu.ru, renat.gumerov\@ksu.ru
\endemail



\vskip 0.5 cm

\abstract
Finite-sheeted covering mappings onto compact connected groups are studied.
It is shown that a finite-sheeted covering mapping \ $p:X\to G$ \ from
a~connected Hausdorff topological space \ $X$ \ onto a compact connected
abelian group \ $G$ \ must be a homeomorphism provided that the character group
of \ $G$ \ admits division by degree of \ $p$.
Using this result, we obtain criteria of triviality for finite-sheeted
covering mappings onto \ $G$ \ in terms of its character group and means on
\ $G$.
In order to establish these facts we show that for a $k$-sheeted covering mapping
$p$ from \ $X$ \ onto a compact connected (in general, non-abelian) group
there exists a~topological group structure on the space $X$ such
that \ $p$ \ becomes a homomorphism of groups. In proving this
result we construct an inverse system of \ $k$-sheeted covering mappings
onto Lie groups which approximates the covering mapping \ $p$.

\endabstract

\endtopmatter

\document

\head
Introduction
\endhead

This paper deals with finite-sheeted covering mappings onto compact
connected groups. Particular attention is paid to covering mappings
onto abelian groups.

The motivation for our work comes partly from the
theory of polynomials and algebraic equations with functional
coefficients. Initially we were interested in pro\-perties
of separable polynomials over the Banach algebra $C(G)$ of all continuous
complex-valued functions on a compact connected abelian group $G$.
A polynomial $ r(w)=w^k + f_1w^{k-1} + f_2w^{k-2} + \dots + f_k$ of degree
\ $k$ \ in the variable \ $w$ \
with coefficients \  $f_1, f_2, \dots, f_k$ \ from \ $C(G)$ \ is said to
be {\it separable } provided that for each element \ $g \in G$ \ the
polynomial \ $w^k + f_1(g)w^{k-1} + f_2(g)w^{k-2} + \dots + f_k(g)$ \
with complex coefficients has no multiple roots in the field of
complex numbers. A separable polynomial $r(w)$ of degree $k$
over the algebra $C(G)$ is said to be {\it totally reducible} if it may be
factored into a product $(w-h_1)\cdot(w-h_2)\cdot \ldots \cdot(w-h_k)$ of
polynomials of degree one, where \ $h_1,h_2,\ldots,h_k \in C(G)$.
Totally reducible polynomials over function algebras and related problems
concerning algebraic equations with functional coefficients have been
studied by E.~A.~Gorin, V.~Ya.~Lin, Yu.~V.~Zyuzin and
others (see, e.g., \cite{1} -- \cite{3}).
In parti\-cular, an immediate consequence of results obtained by
E. A. Gorin and V.~Ya.~Lin (see \cite{1; Corollaries 1.10 -- 1.12}) is
the~following: if the character group \ $\widehat G$ \ is \ $k!$-divisible
then each separable polynomial of degree \ $k$ \ over the algebra \ $C(G)$ \
is totally reducible. One can easily see that the~converse assertion
is also true. It should be noted that the result formulated above can be
interpreted as a generalization of the Walter -- Bohr -- Flanders theorem
concerning algebraic equations with almost periodic coefficients. (\,See
comments in \cite{1; Introduction and \S 1}\,).

As is well known, a separable polynomial \
$r(w)=w^k + f_1w^{k-1} + f_2w^{k-2} + \dots + f_k$ \ over \ $C(G)$ \
generates the $k$-sheeted covering mapping \ $p:X \to G:(g,z) \mapsto g$, \
where the compact space \ $X$ \ is defined as follows:
$$X:=\{(g,z) \in G \times \Bbb C: z^k + f_1(g)z^{k-1}+ \dots + f_k(g)=0\},$$
and the polynomial \ $r(w)$ \ is totally reducible if and only if
the covering mapping $p$ is trivial.

In the present paper it is shown that each $k$-sheeted covering mapping
onto \ $G$ \ is trivial provided that the character group \
$\widehat G$ \ is $k!$-divisible. Thus we have a~criterion for
triviality of all $k$-sheeted covering mappings onto a compact connected
abelian group ( Theorem 3 ).

In studying covering mappings from topological spaces onto
topological groups it is natural to regard the problem on the existence of
a topological group structure on a covering space relative to which
a given covering mapping becomes a~homomorphism of groups. It follows from
classical properties of lifting mappings to covering spaces
that if a covering space is connected and locally path connected, then this
problem has a positive solution, i.e.,
the desired structure exists (\, see, e.g., \cite{4; Theorem 79}\,).
In that case we say that the structure of topological group lifts to
a~covering space, and results of that kind are called the covering group
theorems. In this paper we prove the covering group theorem for
a finite-sheeted covering mapping from a connected Hausdorff topological
space onto a compact connected (\,in general, non-abelian\,) group
(\, Theorem 1\, ). Note that we do not suppose that the group is locally
connected. In order to prove Theorem 1 we first construct a~family of
$k$-sheeted covering mappings onto Lie groups which approximates
a given $k$-sheeted covering mapping (\, see \S 2, Proposition\, ).
Being an object of interest in its own right, Theorem~1 \ serves as
a~tool for studying the structure of covering mappings onto topological
groups (~see \S 4\,).

There are a number of papers treating of covering spaces in which
the~hypothesis of local connectedness is dropped
(~see, e.g., references in \cite{5; Chapter V}). In 1972, R.~Fox,
using  ideas and methods of shape theory, extended results of the classical
covering space theory to overlay mappings of arbitrary connected metric
spaces \cite{6}, \cite{7}. Note that each overlay mapping is a covering
mapping, and the converse implication holds in some important cases.
Namely, it holds if a base space of a covering mapping is locally connected
or a covering mapping is finite-sheeted \cite{7; Theorem 3},
\cite{8; Propositions 2.1, 2.2}. In \cite{9}, \cite{10} and \cite{11} \
Fox's results are generalized to connected topological spaces.

The paper is organized as follows. It is divided into 4 sections.
Section~1 \ contains preliminaries.
In Section~2 the account of the approximate construction
is given.
Section~3 is devoted to the proof of the covering group theorem.
In the last section we apply this theorem to the study of
finite-sheeted covering mappings onto compact connected
abelian groups and to the problem on the existence of means on
topological abelian groups.

The authors are very grateful to Professor S. A. Bogatyi for useful
discussions about the theorems in this paper and for advices and remarks
which improved the~first version of the presentation of the results.
We would like to thank Professors A.~Ya.~Helemskii, A.~S.~Mishchenko
and the members of the Alexandrov seminar and the seminars \
" Algebras in Analysis " \ and \ " Topology and Analysis " \
at Moscow State University for attention to our work and for useful
discussions.
We are very grateful to Professor E.~A.~Gorin for valuable comments.

The results of the present paper were announced in \cite{12} and
\cite{13}.

\head
1. Preliminaries
\endhead

In this section we establish some notation and recall some definitions
and facts that will be needed later on.

Throughout this paper all topological spaces are assumed to be Hausdorff.
\linebreak
A {\it neighborhood} is always an open subset of a topological space.
A {\it mapping} between topological spaces and a {\it homomorphism} between
topological groups mean a~continuous function and a continuous
homomorphism, respectively.

As usual, we denote by $\Bbb N$ the set of all positive integers,
by $\Bbb C$ the space of all complex numbers (equipped with the natural
topology), and by $\Bbb C^m$ the Cartesian product of $m$ copies
of $\Bbb C$, where $m\in \Bbb N$. The set $\{1,2, \ldots, m \}$ is denoted by
$\overline{m}$.

Let \ $X, Y_1, Y_2, \ldots  Y_n$ \ be topological spaces. The {\it diagonal of mappings }
\linebreak
$f_j:X\to Y_j$,
where $j \in \overline{n}$, is defined as follows:
$$f_1 \triangle f_2 \triangle \ldots \triangle f_n:
X \to Y_1 \times Y_2 \times \ldots \times Y_n:
x \mapsto (f_1(x),f_2(x), \ldots, f_n(x) ), $$
where $x\in X$ and $Y_1 \times Y_2 \times \ldots \times Y_n$ is the Cartesian product of
spaces $ Y_1, Y_2, \ldots, Y_n $.


We shall now recall the notions of inverse system and
inverse limit in categories of topological spaces and of topological groups.
For more complete treatments of these notions consult, e.g.,
\cite{14; Chapter VIII, \S 3}, \cite{15; \S 2.5} and \cite{16; Chapter 3}.
We denote by $\Cal {COMP}$ the category of compact spaces and mappings,
and by $\Cal {CGR}$ the category of compact groups and homomorphisms.
Let $\Cal C$ be one of these categories and let
$(\Lambda,\prec)$ be a directed set.
Suppose we are given an object $X_{\lambda}$ from $\Cal C$ for each
$\lambda \in \Lambda$ and
a morphism $\pi_{\lambda}^{\mu}:X_{\mu} \to X_{\lambda}$ from $\Cal C$
for any $\lambda, \mu \in \Lambda$ satisfying $\lambda \prec \mu$.
Moreover, let
$\pi_{\lambda}^{\lambda}:X_{\lambda}\to X_{\lambda}$ be the identity morphism
for each $\lambda \in \Lambda$ and
$\pi_{\lambda}^{\nu} = \pi_{\lambda}^{\mu}{\circ}\pi_{\mu}^{\nu}$
for all $\lambda,\mu,\nu\in\Lambda$ such that $\lambda\prec\mu\prec\nu$.
The collection $\{X_{\lambda},\pi_{\lambda}^{\mu},\Lambda\}$
is called an {\it inverse system} in the category $\Cal C$ over the
{\it index set} $\Lambda$.
The morphisms $\pi_{\lambda}^{\mu}$ are called the {\it bonding morphisms}
of $\{X_{\lambda},\pi_{\lambda}^{\mu},\Lambda\}$. An {\it inverse limit}
$(X,\{\pi_{\lambda}\})$ of the inverse system
$\{X_{\lambda},\pi_{\lambda}^{\mu},\Lambda\}$ consists of an object $X$
from $\Cal C$ and a~family
$\{\pi_{\lambda}:X \to X_{\lambda} \mid \lambda \in \Lambda \}$ of
morphisms from $\Cal C$ having the following {\it universal property}:

1) $\pi_{\lambda}=\pi_{\lambda}^{\mu}{\circ}\pi_\mu$ whenever
$\lambda \prec \mu$;

2) if \ $\tau_{\lambda}: Y \to X_{\lambda}, \quad
\lambda \in \Lambda$, are morphisms from \ $\Cal C$ \ such that \
$\tau_{\lambda}=\pi_{\lambda}^{\mu}{\circ}\tau_\mu$ \ whenever \
$\lambda \prec \mu$, then there exists a unique morphism \ $\rho: Y \to X$ \
from \ $\Cal C$ \ satisfying \ $\tau_{\lambda}=\pi_{\lambda}{\circ}\rho$ \
for each \ $\lambda \in \Lambda$.
\newline
The inverse limit \ $(X,\{\pi_{\lambda}\})$ \ of the inverse system
$\{X_{\lambda},\pi_{\lambda}^{\mu},\Lambda\}$ is denoted by
$(X,\{\pi_{\lambda}\})= \varprojlim\{X_{\lambda},\pi_{\lambda}^{\mu},\Lambda\}$.
One often calls \ $X$ \ itself an inverse system. The~morphisms \
$\pi_\lambda$ \ are called the {\it projections}.
In the category \ $\Cal C$ \ an inverse limit \
$\varprojlim\{X_{\lambda},\pi_{\lambda}^{\mu},\Lambda\}$ \
exists. Namely, the object \ $X_\infty$ \ from \ $\Cal C$ \ consisting of all
threads of the inverse system \ $\{X_{\lambda},\pi_{\lambda}^{\mu},\Lambda\}$ \
and the family of the canonical projections of \ $X_\infty$ \ to \ $X_\lambda$ \
constitute the inverse limit. If there are two inverse limits of
the same inverse system in some category, then they are isomorphic:
the morphisms from the~universal properties for these inverse limits are
isomorphisms \ \cite{16; Corollary 3.2}.

A {\it morphism of inverse systems} \
$\{X_{\lambda},\pi_{\lambda}^{\mu},\Lambda\}$ \ and \
$\{Y_{\lambda},\sigma_{\lambda}^{\mu},\Lambda\}$ \ in \ $\Cal C$ \
is a family
$\{\tau_\lambda: X_\lambda \to Y_\lambda \mid \lambda \in \Lambda\}$
of morphisms from $\Cal C$ such that $\tau_\lambda \circ \pi_{\lambda}^{\mu} =
\sigma_{\lambda}^{\mu} \circ \tau_\mu$ whenever $\lambda \prec \mu$.
It induces a morphism $\tau_{\infty}:X_{\infty}\to Y_\infty$ in $\Cal C$
given by $\tau_{\infty}(\{x_{\lambda}\}_{\lambda\in\Lambda}) =
\{\tau_{\lambda}(x_{\lambda})\}_{\lambda\in\Lambda}$,
where $\{x_{\lambda}\}_{\lambda\in\Lambda}\in X_\infty$.
The morphism $\tau_{\infty}$ is called the {\it limit morphism induced by}
$\{\tau_\lambda: \lambda \in \Lambda\}$. We say that a morphism
$\tau:X \to Y$ from $\Cal C$ is {\it up to isomorphism}
the limit morphism induced by the morphism of inverse systems
$\{\tau_\lambda: \lambda \in \Lambda\}$ provided that there are isomorphisms
$\rho: X \to X_\infty$ and $\sigma: Y \to Y_\infty$ from $\Cal C$
such that $\tau_{\infty} \circ \rho=\sigma \circ \tau$.

Recall that a surjective mapping $p:X\to Y$ between topological spaces $X$
and $Y$ is called a {\it finite-sheeted covering mapping} \ if \ it
is a $k$-{\it sheeted} \ (\, $k$-{\it fold}\,) {\it covering mapping}
for some $k\in \Bbb N$.
That is, for every point $y\in Y$, there exists a~neighborhood $W$ in $Y$
and a partition of the inverse image $p^{-1}(W)$ into neighborhoods
$V_1,V_2,\ldots,V_k$ in $X$ such that, for each $n\in \overline{k}$,
the restriction of $p$ to $V_n$ is a homeomorphism of $V_n$ onto $W$.
The integer $k$ is called a {\it degree} of the mapping $p$.
The~neighborhood $W$ is said to be {\it evenly covered} \ by \ $p$,
and the collection $\{V_n: n\in \overline{k} \}$ is called a partition of
$p^{-1}(W)$ into {\it slices}.
(\, As usual, a partition of an arbitrary set $A$ is a collection
of disjoint subsets of $A$ whose union is all of $A$.)

A $k$-sheeted covering mapping \ $p:X\to Y$ \ is said to be {\it trivial}
if there exists a~homeomorphism \ $\rho$ \ from the Cartesian product \
$Y\times \overline{k}$, where \ $\overline{k}$ \ is equipped with
the discrete topology, to \ $X$ \
such that \ $p \circ \rho(y,n)=y$ \ for any \ $(y,n)\in Y\times \overline{k}$.

Suppose we are given a $k$-fold covering mapping \ $p:X\to Y$ \ from
a non-connected space \ $X$ \ onto a connected space \ $Y$, where $k\geq 2.$
It is a well-known fact that for any component \ $X_0$ \ of the space \ $X$ \
the restriction \ $p\vert_{X_0}:X_0\to Y$ \ of the~mapping \ $p$ \ to the
subspace \ $X_0$ \ is a finite-sheeted covering mapping whose degree is less
than \ $k$.

Throughout the paper $G$ stands for a compact connected group,
and $p:X\to G $ is a $k$-fold covering mapping, where $X$ is a topological
space. Note that the space $X$ is necessarily compact.

Let us consider an inverse system
$\{\,G_{\lambda},\pi_{\lambda}^{\mu},\Lambda\,\}$ in $\Cal {CGR}$
such that $(G, \{\pi_{\lambda}\})=
\varprojlim \{G_{\lambda},\pi_{\lambda}^{\mu},\Lambda\}$,
where all $G_\lambda$ are connected Lie groups, and all bonding
homomorphisms and projections are open and surjective.
Moreover, if $G$ is abelian, then each $G_{\lambda}$ is assumed to be
abelian as well (see [17; \S25]).

A family of all sets $\pi_{\lambda}^{-1}(U)$, where
$U$ is a neighborhood in $G_{\lambda}$ and $\lambda$ runs over a~subset
cofinal in $\Lambda$, is a basis for a topology on $G$
\cite{15; Proposition 2.5.5}. For $\nu\in\Lambda$ we denote by $\Lambda_\nu$
the set $\{\,\lambda\in\Lambda: \nu\prec\lambda\,\}.$
Then
$(G,\{\pi_\lambda\}_{\lambda \in \Lambda_\nu})=\varprojlim{\{\,G_\lambda,\pi_{\lambda}^{\mu},\Lambda_\nu\,\}}$
\cite{14; Chapter VIII, \S 3}.

Using compactness of the group $G$ and the above-mentioned description of
\
a~basis for the topology on $G$, one can easily see that there is a finite open
cover $\{\,W_n: n \in \overline{m}\,\}$ of $G$ (i.e., a finite family of
neighborhoods in $G$ whose union is the~whole group $G$) satisfying
the following conditions. Every neighborhood $W_n$ is an~evenly covered
by the covering mapping $p$,  and there exists an index $\alpha\in\Lambda$
and a~finite open cover $\{\,W_n^{\alpha}: n \in \overline{m}\,\}$ of
$G_\alpha$ such that $W_n=\pi_{\alpha}^{-1}(W_{n}^{\alpha})$ for each
$n \in \overline{m}$. From now on we fix $\alpha\in\Lambda$ and the indexed
family $W_1,W_2,\dots,W_m$. In addition, for each set $p^{-1}(W_n)$
we fix its indexed partition into slices
$V_{1}^{n},V_{2}^{n},\dots,V_{k}^{n}$. Thus for each $n \in \overline{m}$
we have the disjoint union
$$
p^{-1}(W_n)=\bigsqcup_{l=1}^{k} V_{l}^{n} \tag {$1$}
$$
such that for each $l \in \overline{k}$ the neighborhood $V_{l}^{n}$ is
mapped homeomorphically onto the neighborhood $W_{n}$ under $p$.

\head
2. Approximation of covering mapping
\endhead


In the present section we shall construct a family of \ $k$-sheeted covering
mappings onto \ $G_\lambda$'s \ which approximates the covering mapping
\ $p:X\to G$. Roughly speaking, it turns out that in an appropriate
category of arrows the covering mapping \ $p:X\to G$ \ is an inverse
limit of an inverse system consisting of $k$-sheeted covering
mappings onto $G_\lambda$'s.

Throughout this section we "forget"
about group structures of objects,
so that $G$, $G_\lambda$ and $\pi_\lambda^\mu, \pi_\lambda$ are
considered as objects and morphisms of the category $\Cal {COMP}$,
i.e., as only compact spaces and mappings respectively. (\, In other words,
we apply the forgetful functor from $\Cal {CGR}$ to $\Cal {COMP}$\,).

We need first to define functions and neighborhoods which together
with those we have introduced in the preceding section will play the crucial
role in our construction.

It is clear that one can "shrink" the cover $\{\,W_n: n \in \overline{m}\,\}$
to a finite open cover $\{\,U_n: n \in \overline{m}\,\}$ of the space $G$
such that
$$
U_n = \pi_{\alpha}^{-1}(U^{\alpha}_n) \qquad \text{and} \qquad \overline U_n \subset W_n
\qquad \text{for each} \qquad n \in \overline{m}, \tag {$2$}
$$
where $\{\,U_{n}^{\alpha}: n \in \overline{m}\,\}$ is a finite open cover
of the space $G_\alpha$ which refines the cover
$\{\,W_{n}^{\alpha}: n \in \overline{m}\,\}$ and satisfies the condition
$U_{n}^{\alpha} \subset W_{n}^{\alpha}$ \ for each $n$.

Let $\{\,\phi_n: n \in \overline{m} \,\}$ be a partition of unity on
the space $G_\alpha$ dominated by the~cover
$\{\,U_{n}^{\alpha}: n \in \overline{m}\,\}$, that is, an indexed family
of mappings
$$
\phi_n: G_\alpha\to [0,1] \quad \text{for} \quad n \in \overline{m}
$$
having the following properties:

---  the support of $\phi_n$, i.e., the closure of the set \
$\phi_n^{-1}((0,1])$, is contained in the~neighborhood \ $U_{n}^{\alpha}$ \
for each $n$;

---  the equality \ $\sum_{n=1}^{m}\phi_n(g) =1$ \ holds for each
$g \in G_\alpha$.

For every $n \in \overline{m}$ \ we consider the compositions of mappings
defined as follows:
$$\psi_n := \phi_n\circ\pi_\alpha : G\to [0,1] \quad \text{and} \quad
\widehat\psi_n := \phi_n\circ\pi_\alpha{\circ}p : X\to [0,1].$$

It is straightforward to check that the families of mappings
$\{\,\psi_n: n \in \overline{m}\,\}$ \ and \
$\{\,\widehat\psi_n: n \in \overline{m}\,\}$ are the partitions of unity
on $G$ and on $X$ dominated by the covers
$\{\,U_n: n \in \overline{m}\,\}$ \ and \  $\{\,p^{-1}(U_n): n \in \overline{m}\,\}$,
respectively.

We next define functions $f_n : X\to\Bbb C$ \ and \
$\widehat f_n : X\to\Bbb C$ for each $n\in \overline{m}$ by setting
$$
f_n(x) :=\left\{
\aligned  0,    \quad \text{if} \quad x\notin p^{-1}(W_n)=\bigsqcup_{l=1}^k V_{l}^n; \qquad \qquad \\
\exp(i\frac{2\pi}{k}(l-1)), \quad \text{if} \quad x\in V_{l}^n, \quad l \in
\overline{k}; \qquad (i^2 = -1);
\endaligned \right. \tag {$3$}
$$

$$
\widehat f_n(x):=\widehat \psi_n(x)f_n(x)  \qquad \text{for} \quad x\in X.
$$
In view of (2), the functions $\widehat f_1,\widehat f_2,\dots,\widehat f_m$
are continuous at each point of $X.$

We claim that the diagonal of mappings
$$p\triangle\widehat f_1\triangle\widehat f_2\triangle\cdots\triangle\widehat f_m :
X\to G\times\Bbb C^m$$
is a homeomorphic embedding, i.e., the space $X$ is mapped homeomorphically
onto the image space \
$p\triangle\widehat f_1\triangle\widehat f_2\triangle\cdots\triangle\widehat f_m(X)
\subset G\times\Bbb C^m$ \
under the~diagonal of the~mappings $p,\widehat f_1,\ldots,\widehat f_m$.
To prove this it suffices to show that
$p\triangle\widehat f_1\triangle\widehat f_2\triangle\cdots\triangle\widehat f_m$
is injective (because it is continuous and the space $X$
is compact), or, equivalently, that the family of mappings
$\{p,\widehat f_1,\widehat f_2,\dots,\widehat f_m\}$ separates points.
For this purpose we take two distinct points $x, y \in X$.
If $p(x)=p(y)$, then $\widehat \psi_n(x)=\widehat \psi_n(y)$ for each
$n\in \overline{m}$. Since the functions $\widehat \psi_n$ form a partition
of unity on $X$, we have
$$
\sum\limits_{n=1}^m \widehat \psi_n(x)=1,
$$
which implies $\widehat \psi_n(x)>0$ for some $n$. Hence the points $x$
and $y$ belong to the set $p^{-1}(U_n)\subset p^{-1}(W_n)$.
But these points lie in different slices of $p^{-1}(W_n)$, so that
$f_n(x)\ne f_n(y)$. Therefore, we have $\widehat f_n(x)\ne\widehat f_n(y)$
as well. Thus the diagonal of mappings
$p\triangle\widehat f_1\triangle\widehat f_2\triangle\cdots\triangle\widehat f_m$
is an injection. Consequently, it is a homeomorphic embedding, as
claimed.

Now for each index $\lambda\in\Lambda_\alpha$ let us define a compact
space $X_\lambda$ by letting
$$X_\lambda :=\bigl\{\,\bigl(\pi_\lambda(p(x)),\widehat f_1(x),\widehat f_2(x),\cdots,\widehat f_m(x)\bigr)
: x\in X\,\bigr\} \subset G_\lambda\times\Bbb C^m,$$
i.e., $X_\lambda$ is the image of the space $X$ under the following diagonal
of mappings
$$(\pi_\lambda{\circ}p)\triangle
\widehat f_1\triangle\widehat f_2\triangle\cdots\triangle\widehat f_m:
X \to G_\lambda\times\Bbb C^m.$$
We denote by $h_\lambda: X\to X_\lambda$ the surjective mapping given by
$$
h_\lambda (x)=(\pi_\lambda{\circ}p)\triangle
\widehat f_1\triangle\widehat f_2\triangle\cdots\triangle\widehat f_m (x)
\qquad \text{for} \quad x\in X,$$
and by $p_\lambda$  the projection of $X_\lambda$ onto the first
coordinate, that is,
$$p_\lambda: X_\lambda \to G_\lambda:
(\pi_\lambda(p(x)),\widehat f_1(x),\widehat f_2(x),\cdots,\widehat f_m(x))
\mapsto \pi_\lambda(p(x)), \quad x\in X.$$

Obviously, for each index $\lambda \in \Lambda_\alpha$, the following diagram
$$
\CD
X   @>p>>   G \\
@Vh_{\lambda}VV   @VV\pi_{\lambda}V \\
X_\lambda @>p_{\lambda}>> G_{\lambda}
\endCD \tag{$4$}
$$
is commutative, i.e., $\pi_{\lambda}\circ p=p_{\lambda}\circ h_\lambda$.
In particular, as an immediate consequence of this fact,
the composition $p_{\lambda}\circ h_\lambda$ and
the mapping $p_\lambda$ are open.

Let us fix an index $\lambda \in \Lambda_\alpha$. For an arbitrary
$g\in G_\lambda$, commutativity of the~diagram (4) yields the equality
$$
h_\lambda^{-1}(p_\lambda^{-1}(g))=(\pi_\lambda\circ p)^{-1}(g),
$$
which implies
$$
p_\lambda^{-1}(g)=\{(g,\widehat f_1(x),\widehat f_2(x),\dots,\widehat f_m(x)):x\in
(\pi_\lambda\circ p)^{-1}(g)\}. \tag{$5$}
$$
Furthermore, for each $n \in \overline{m}$, we have:
$$
\widehat\psi_n(x)=\phi_n\circ\pi_{\alpha}^{\lambda}(g) \quad
\text{whenever} \quad x\in (\pi_\lambda{\circ}p)^{-1}(g).
$$
In other words, the mapping $\widehat\psi_k(x)$ is constant on the fiber
$(\pi_\lambda{\circ}p)^{-1}(g).$ Therefore we can rewrite (5) as follows:
$$
p_{\lambda}^{-1}(g)=\{(g,\phi_1(\pi_{\alpha}^{\lambda}(g))f_1(x),\ldots,\phi_m(\pi_{\alpha}^{\lambda}(g))f_m(x)) :
x\in (\pi_\lambda{\circ}p)^{-1}(g)\}. \tag {$6$}
$$

But, for each $n\in \overline{m}$, the cardinal number of the set
$\{f_n(x): x\in (\pi_\lambda{\circ}p)^{-1}(g)\}$
is finite. Namely, it is at most $k+1$ (\,see(3)\,).
Using this fact and the equality (6), we conclude that the set
$p_{\lambda}^{-1}(g)$  is finite.
Generally speaking, the cardinal number of the fiber $p_{\lambda}^{-1}(g)$
depends on $g$. So we can not assert that the mapping $p_{\lambda}$ is
a~finite-sheeted covering.

Nevertheless, it turns out that the following lemma is valid.

\proclaim{Lemma}
There exists an index $\beta\in\Lambda_\alpha$ such that the mapping
$p_\lambda : X_\lambda\to G_\lambda$ is a~$k-fold$ covering mapping for
each $\lambda\in\Lambda_\beta$.
\endproclaim
\demo{Proof} The proof consists of three steps. The first step is to
construct a finite open cover $ \{\,O_s: s \in \overline{t}\,\}, \
t\in \Bbb N$, \ of the space $G$ which refines the cover
$\{\,W_n: n \in \overline{m}\,\}$ and satisfies some additional
requirements. We choose $\beta\in\Lambda_\alpha$ during this construction.
Then we use the cover $ \{\,O_s: s \in \overline{t}\,\}$ to show that
the mappings $p_\lambda,\lambda\in\Lambda_\beta$, are $k$-fold
covering mappings.

{\it Step 1.} \ \  We prove the following: \
{\it \ There exists a finite open cover \ \ $ \{\,O_s: s \in \overline{t}\,\}$
of the space $G$ which has the following properties:}

$O1)$ {\it For each $s \in \overline{t}$ we have
$O_s = \pi_{\beta}^{-1}(O_{s}^{\beta}),$
where $\beta\in\Lambda_\alpha,$
and the family $\{\,O_{s}^{\beta}:s \in \overline{t}\,\}$ is an open cover
of the space $G_\beta;$}

$O2)$ {\it For each $s \in \overline{t}$ there exists a partition
of the set of positive integers $\overline{m}$:
$$
\overline{m} = \{\,a_1,a_2,\ldots,a_r\,\} \ \sqcup \ \{\,b_1,b_2,\dots,b_{m-r}\,\},
\tag {$7$}
$$
such that $$O_s\subset\bigcap_{j=1}^{r}W_{a_j} \quad \text{and} \quad
O_s\bigcap(\bigcup_{j=1}^{m-r}\overline U_{b_j}) = \emptyset.
$$
Furthermore, all mappings \ $\widehat\psi_{b_j}$
vanish on the set \ $p^{-1}(O_s)$,
i.e., \ $\widehat\psi_{b_j}(x)=0 $ whenever $ x\in p^{-1}(O_s)$
and $j \in \overline{m-r},$ \ and the image sets
$f_{a_1}{\triangle}f_{a_2}{\triangle}\cdots{\triangle}f_{a_r}(p^{-1}(g))$
consist of the same $k$ elements while $g$ runs over $O_s.$}

(\, To simplify the notation, we make use of the symbols $a_j,b_j$ and
$r$ instead of $a_{sj},b_{sj}$ and $r_{s}$ respectively. Of course, if $r=m,$
then the set \ $\{\,b_j: j \in \overline{m-r}\,\}$ \ is assumed to
be the empty set $\emptyset$\,).

Taking an arbitrary element $g\in G$, we consider the partition
$$\{\,W_n: n \in \overline{m}\,\}=
\{\,W_{a_j}: j \in \overline{r}\,\} \, \bigsqcup \, \{\,W_{b_j}:j \in
\overline{m-r}\,\},$$
which is uniquely determined by the following requirements:
$$g\in\bigcap_{j=1}^{r}W_{a_j} \qquad
\text{and} \qquad g\notin\bigcup_{j=1}^{m-r}W_{b_j}.$$
In other words, we have the partition (7) of the set of positive integers
$\overline{m}$.

It follows immediately from (2) that
$$g\notin\bigcup_{j=1}^{m-r}\overline U_{b_j}. \tag {$8$}$$

Let $x$ be an element of the set $ p^{-1}(g)$. For each integer
$a_j$ from (7), there is a~unique slice
in the partition (1) of $p^{-1}(W_{a_j})$ into slices which contains $x.$
A~subscript of that slice is denoted by $l_j$. Thus the condition $x\in
V_{l_j}^{a_j}$ is fulfilled.

We next define a neighborhood $O(x)$ of the point $x$ by setting
$$O(x):=\bigcap_{j=1}^{r}V_{l_j}^{a_j}.$$
In view of (3), for each $j\in \overline{r}$, we have:
$$f_{a_j}(y)=\exp(i\frac{2\pi}{k}(l_j - 1))\quad
\text{whenever} \quad y\in O(x).$$

In a similar way, assuming that $p^{-1}(g)=\{\,x_1,x_2,\dots,x_k\,\}$,
we define disjoint neighborhoods \ $O(x_1),O(x_2),\dots,O(x_k)$ \ of
the points \ $x_1,x_2,\cdots,x_k$ \ respectively. It is clear that for each
$l \in \overline{k}$ the restriction of the diagonal of
mappings $f_{a_1}{\triangle}f_{a_2}{\triangle}\cdots{\triangle}f_{a_r}$
to the neighborhood $O(x_l)$ is constant.

Using (8), we choose a neighborhood $O(g)$ of the point $g$ such that
$$O(g)\subset\bigcap_{l=1}^{k}p(O(x_l)) \qquad \text{and} \qquad
O(g)\bigcap(\bigcup_{j=1}^{m-r}\overline U_{b_j})=\emptyset.$$
Obviously, we may assume without loss of generality that
$O(g)=\pi_{\lambda}^{-1}(U),$ where
$\lambda\in\Lambda_\alpha$, and $U$ is a nonempty neighborhood in
$G_\lambda.$

One can easily see that all functions \
$\widehat\psi_{b_1},\widehat\psi_{b_2},\dots,\widehat\psi_{b_{m-r}}$ \
vanish on the set $p^{-1}(O(g))$ \ (\, provided that the set \
$\{\,b_1,b_2,\dots,b_{m-r}\,\}$ \ is nonempty \,), and the~image set
$f_{a_1}{\triangle}f_{a_2}{\triangle}\cdots{\triangle}f_{a_r}(p^{-1}(O(g)))$
consists of $k$ points. Moreover, for any $g^{\prime}\in O(g)$ with
$p^{-1}(g^{\prime})=\{y_1,y_2,\dots,y_k\}$ the following equality holds:
$$f_{a_1}{\triangle}f_{a_2}{\triangle}\cdots{\triangle}f_{a_r}(p^{-1}(O(g)))=
\{\,f_{a_1}{\triangle}f_{a_2}{\triangle}\cdots{\triangle}f_{a_r}(y_l):
l \in \overline{k}\,\}.$$

Since the space $G$ is compact, there exists a finite open cover
$$
\bigl\{\,O(g_s)=\pi_{\lambda_s}^{-1}(U_{s}^{\lambda_s}): g_s \in G, U_{s}^{\lambda_s}
\quad \text{is a neighborhood in} \quad G_{\lambda_s},\lambda_s \in
\Lambda_\alpha, s \in \overline{t}\,\bigr\}
$$
of \ $G$ \ satisfying the condition \ $O2).$ The neighborhood \ $O(g_s)$ \ is
denoted by $O_s$.

Now we choose an index $\beta\in\Lambda_\alpha$ such that
$\lambda_s\prec\beta$ for all $s \in \overline{t}.$ It is clear that
$$O_s=\pi_{\beta}^{-1}(O_{s}^{\beta}), \quad \text{where} \quad
O_{s}^{\beta}=(\pi_{\lambda_s}^{\beta})^{-1}(U_{s}^{\lambda_s}), \quad
s \in \overline{t}.$$

Since the mapping $\pi_\beta$ is surjective, the family
$\{\,O_{s}^{\beta}: s \in \overline{t}\,\}$ is an open cover of $G_\beta.$
Thus the cover $\{\,O_s: s \in \overline{t}\,\}$ also has the property $O1)$,
as desired.

{\it Step 2.} Now we show that the open mapping \
$p_\beta: X_\beta\to G_\beta$ \ is a $k$-fold covering mapping.
Clearly, to prove this fact it is enough to show that the mapping \ $p_\beta$ \
is \ $k$-to-$1$, i.e., the fiber \ $p_{\beta}^{-1}(g)$ \ consists of \ $k$ \
points for each \ $g\in G_\beta$.

Choose any \ $g\in G_\beta.$ For the set \ $p_{\beta}^{-1}(g)$ \ we have
the equalities (5) and (6) (with $\lambda=\beta$).

Since the family $\{\,O_{s}^{\beta}: s \in \overline{t}\,\}$
constructed in Step 1 is a cover of the space $G_\beta$, there is
an index $s$ such that $\pi_{\beta}^{-1}(g)\subset O_s.$  We fix such $s$
and consider the partition (7) corresponding to this index.

If the set $\{\,b_j:j \in n_{m-r}\,\}$ is nonempty, then
(see the property $O2)$ in Step 1) for each $b_j$ we have:
$$\widehat f_{b_j}(x)=\widehat \psi_{b_j}(x)f_{b_j}(x)=0 \quad \text{whenever}
\quad x\in (\pi_\beta{\circ}p)^{-1}(g).$$
Hence the sets $$p_{\beta}^{-1}(g) \quad \text{and} \quad
\bigl\{\,\widehat f_{a_1}\triangle\widehat f_{a_2}\triangle\cdots\triangle\widehat f_{a_r}(x):
x\in (\pi_\beta{\circ}p)^{-1}(g)\,\bigr\}$$ have the same cardinality.

Now we fix a point  $x_0 \in (\pi_\beta{\circ}p)^{-1}(g) $. As was noted,
for each integer $a_j$ from the partition \ (7) \ corresponding to $s$, we have:
$$\widehat\psi_{a_j}(x)=\widehat\psi_{a_j}(x_0) \quad \text{whenever} \quad
x\in(\pi_\beta{\circ}p)^{-1}(g).$$
Therefore, it suffices to show that the set
$$
\bigl\{\,(\widehat\psi_{a_1}(x_0)f_{a_1}(x),\widehat\psi_{a_2}(x_0)f_{a_2}(x),\dots,
\widehat\psi_{a_r}(x_0)f_{a_r}(x)) : x\in (\pi_\beta{\circ}p)^{-1}(g)\,\bigr\} \tag {$9$}
$$
consists of \ $k$\ \ points.

The family $\{\,\widehat\psi_n: n \in \overline{m}\,\}$ is the partition of
unity on the space $X$, and the~mappings \
$\widehat\psi_{b_1},\widehat\psi_{b_2},\dots,\widehat\psi_{b_{m-r}}$ \
vanish on the set $(\pi_\beta{\circ}p)^{-1}(g)$, whence
$$\sum_{j=1}^{r}\widehat\psi_{a_j}(x_0) =1,$$
which implies the inequality
$\widehat\psi_{a_j}(x_0)>0$ for at least one $a_j.$

Using the last observation together with the property $O2)$ from Step 1 and
the~fact that each function $f_{a_j}$ separates points of every fiber
$p^{-1}(g^{\prime})$, where $g^{\prime}\in \pi_{\beta}^{-1}(g)$, one can
easily see that the set (9) consists of $k$ points. Thus the same is valid
for the set $p_{\beta}^{-1}(g)$.

{\it Step 3}. To show that all mappings $p_\lambda : X_\lambda\to G_\lambda$,
where $\lambda\in\Lambda_\beta$, are $k$-fold covering mappings we first fix
an index $\lambda\in\Lambda_\beta.$

Then we define the family $\{\,O_{s}^{\lambda}: s \in \overline{t}\,\}$ of
neighborhoods by letting
$$O_{s}^{\lambda} :=(\pi_{\beta}^{\lambda})^{-1}(O_{s}^{\beta}),$$
where \ $\{\,O_{s}^{\beta}: s \in \overline{t}\,\}$ \ is the open cover of
the space $G_\beta$ from $O1)$.
Obviously, the~family \ $\{\,O_{s}^{\lambda}: s \in \overline{t}\,\}$ \
is an open cover of the space $G_\lambda$. In addition, the~equality \
$\pi_\beta=\pi_\beta^{\lambda}\circ \pi_\lambda$ \ implies \
$O_s =\pi_{\lambda}^{-1}(O_{s}^{\lambda})$ \
for each $s \in \overline{t}$.

Applying Step 2 to the mapping $p_\lambda : X_\lambda\to G_\lambda$ one
obtains the desired conclusion. This completes the proof of the lemma.
\qed
\enddemo

Now for any $\lambda,\mu\in\Lambda_\beta$
satisfying $\lambda\prec\mu,$ we define the mapping
$h_{\lambda}^{\mu} : X_\mu\to X_\lambda$ by the rule
$$h_{\lambda}^{\mu}(h_\mu(x))  = h_\lambda(x),\quad \text{for} \quad x\in X.$$
Clearly, this mapping is well-defined.

One checks at once that the collection \
$\{\,X_\lambda,h_{\lambda}^{\mu},\Lambda_\beta\,\}$ \
is an inverse system in the category $\Cal {COMP}$, and
$\{\,p_\lambda:\lambda\in\Lambda_\beta\,\}$ is a morphism of the inverse
systems \ $\{\,X_\lambda,h_{\lambda}^{\mu},\Lambda_\beta\,\}$ \  and
\ $\{\,G_\lambda,\pi_{\lambda}^{\mu},\Lambda_\beta\,\}$. Thus each diagram
$$
\CD
X_\lambda    @<h_{\lambda}^{\mu}<<  X_\mu \\
@Vp_{\lambda}VV                  @VVp_{\mu}V \\
G_\lambda     @<\pi_{\lambda}^{\mu}<<  G_\mu,
\endCD \tag {$10$}
$$
is commutative.

Let us consider the limit morphism \ $p_\infty:X_\infty \to G_\infty$ \
induced by \ $\{p_\lambda:\lambda\in\Lambda_\beta\}$. Denote by
$\sigma : G\to G_\infty$ the homeomorphism from the universal property
for  $G_\infty$ and by $\rho: X\to X_\infty$ the mapping
from the universal property for  $X_\infty$.
It is clear that $\sigma(g) = \{\pi_\lambda(g)\}_{\lambda\in\Lambda_\beta},$
where $g\in G$, and $\rho(x) = \{h_\lambda(x)\}_{\lambda\in\Lambda_\beta}$,
where $x \in X$.

Commutativity of the diagram (4) yields commutativity of the following
diagram:
$$
\CD
X   @>\rho>>    X_{\infty} \\
@VpVV              @VVp_{\infty}V \\
G   @>\sigma>>  G_{\infty}
\endCD
\tag {$11$}
$$

Moreover, it is easy to check that $\rho$ is a homeomorphism.
For the sake of completeness, we show this below.

Let $x,y\in X$ and $\rho(x) = \rho(y).$ Then
$\widehat f_j(x) = \widehat f_j(y)$ for each $j\in \overline{m}$ and
$\pi_\lambda(p(x)) = \pi_\lambda(p(y))$ for each $\lambda\in\Lambda_\beta.$
Since the family of projections $\{\pi_\lambda;\lambda\in\Lambda_\beta\}$
separates points of the space $G$, we obtain the equality
$p(x) = p(y)$. Therefore
$$p\triangle\widehat f_1\triangle\widehat f_2\triangle\ldots\triangle\widehat f_m(x)=
p\triangle\widehat f_1\triangle\widehat f_2\triangle\ldots\triangle\widehat f_m(y).$$
But the diagonal \
$p\triangle\widehat f_1\triangle\widehat f_2\triangle\ldots\triangle\widehat f_m$ \
is injective, so that $x=y$ and $\rho$ is injective as well.

Now we prove that $\rho :X\to X_\infty$ is surjective. For fixed
$x=\{x_\lambda\}_{\lambda\in \Lambda_\beta}\in X_\infty,$ where
$x_\lambda\in X_\lambda$ for each $\lambda\in \Lambda_\beta$, we show that
$x = \rho(y)$ for some point $y\in p^{-1}(g),$ where $g$ is the unique point
of the space $G$ such that $\sigma(g) = p_\infty (x)$.

For each $\lambda\in \Lambda_\beta$ we choose a point \ $z_\lambda\in X$ \
such that $x_\lambda=h_\lambda(z_\lambda)$.
We have the equality
$h_\beta(z_\beta) = h_{\beta}^{\lambda}(h_\lambda(z_\lambda))$ which
implies that
$$
h_\lambda(z_\lambda) = (\pi_\lambda(p(z_\lambda)),\widehat f_1(z_\beta),\widehat f_2(z_\beta),
\ldots,\widehat f_m(z_\beta)). \tag {$12$}
$$
Thus the last \ $m$ \ coordinates of all \ $h_\lambda(z_\lambda)$ \ are the same.

Let \ $p^{-1}(g) = \{y_1,y_2,\ldots,y_k\}.$
One can readily verify that
$$
(\pi_\beta{\circ}p)(y_n) = (\pi_\beta{\circ}p)(z_\beta)
\quad \text{for each} \quad y_n\in p^{-1}(g), \quad n \in \overline{m}.
\tag {$13$}
$$

This implies that \ $z_\beta,y_1,y_2,\ldots,y_k\in p^{-1}(O_s)$ \ for
some neighborhood $O_s$ from the cover \ $\{O_1,O_2,\ldots,O_t\}$ \ of the
space \ $G$ \ constructed in the proof of Lemma (\,see Step 1\,).
We fix such $O_s$ and consider the partition (7) of the set $\overline{m}$
corresponding to this neighborhood. By the property $O2)$, there is a unique
point $y_l$ in the set $p^{-1}(g)$ satisfying the equality
$$
f_{a_1}{\triangle}f_{a_2}{\triangle}\ldots{\triangle}f_{a_r}(y_l)
= f_{a_1}{\triangle}f_{a_2}{\triangle}\ldots{\triangle}f_{a_r}(z_\beta).
\tag {$14$}
$$

We claim that \ $\rho(y_l)=x$.

Making use of the equalities (12), (13) and the property $O2)$, we obtain:
$$
\widehat f_{a_j}(y_l) = \widehat f_{a_j}(z_\beta)
\quad \text{for} \quad j\in \overline{r}
$$
and, provided that the set $\{b_j:j\in \overline{m-r}\}$ is nonempty,
$$\widehat f_{b_j}(y_l)  = \widehat f_{b_j}(z_\beta) =0
\quad \text{for} \quad j\in \overline{m-r}.$$

As an immediate consequence of the above observations, we get:
$$
\widehat f_1\triangle\widehat f_2\triangle\ldots\triangle\widehat f_m(y_l)
=\widehat f_1\triangle\widehat f_2\triangle\ldots\triangle\widehat f_m(z_\beta).
\tag 15
$$

One can easily see that
$$
 \pi_{\lambda} {\circ} p(y_l) = \pi_\lambda{\circ}p(z_\lambda)
 \quad \text{for each} \quad \lambda \in \Lambda_\beta.
\tag 16
$$

Finally, combining (12),(15) and (14), we get $\rho(y_l) = x$, as claimed.

Thus $\rho$ is a bijection, so that, being continuous and defined on
a compact space, it is a homeomorphism.

We summarize the results of this section in the following proposition.

\proclaim{Proposition} The covering mapping $p:X\to G$  is up to isomorphism
the limit morphism induced by the morpism 
$$\{p_\lambda: \lambda\in\Lambda_\beta\}:
\{X_\lambda,h_{\lambda}^{\mu},\Lambda_\beta\} \longrightarrow
\{G_\lambda,\pi_{\lambda}^{\mu},\Lambda_\beta\}$$
of inverse systems in the category $\Cal {COMP}$, where $p_\lambda$ is a $k$-fold
covering mapping for each $\lambda \in \Lambda_\beta$.
\endproclaim

\head
3. Covering group theorem
\endhead

In this section we shall show that the structure of topological group on \
$G$ \ lifts to the covering space \ $X$. We begin by recalling the covering
group theorem for groups which are connected and locally path connected
(\,see, e.g., \cite{4, Theorem 79}, \cite{5, Chapter V, Exercises 5.2, 5.4}\,).

\proclaim{Theorem}  \
Let \ $\omega: \widetilde{\Gamma} \to \Gamma $ \ be a covering mapping from
a path connected space \ $\widetilde{\Gamma}$ \ onto a connected locally
path connected topological group $\Gamma$ with identity \ \ $e$.
Then for any point $\widetilde{e} \in \omega^{-1}(e)$ there exists a unique
structure of topological group on \ $\widetilde{\Gamma}$ such that
 \ $\widetilde{e}$ is the identity and \
$\omega: \widetilde{\Gamma} \to \Gamma $ \ is a homomorphism of
topological groups. Furthermore, if \ $\Gamma$ is abelian, then \
$\omega$ \ is a homomorphism of abelian groups.
\endproclaim

Note that the covering mapping $\omega: \widetilde{\Gamma} \to \Gamma $ is not
assumed to be finite-sheeted, and the topological group $\Gamma $ is not
assumed to be compact in the hypothesis of Theorem.

Using Proposition from the previous section and this Theorem, we have:

\proclaim{Theorem 1}
Let \ $p:X \to G $ \ be a finite-sheeted covering mapping from a connected
space \ $X$ \ onto a compact topological group \ $G$ \ with identity \ $e$.
Then for any point \ $\widetilde{e} \in p^{-1}(e)$  \
there exists a unique structure of topological group on \ $X$ \ such that
\ $\widetilde{e}$ \ is the identity and \ $p:X \to G $ \ is a homomorphism
of compact groups. Furthermore, if \ $G$ \ is abelian, then \
$p$ \ is a homomorphism of abelian groups.
\endproclaim

\demo{Proof}
Let \ $\{\,p_{\lambda}: \lambda\in\Lambda_\beta\,\}$ \ be the morphism of
inverse systems in the category \ $\Cal {COMP}$ \ constructed in the previous
section. We now turn it into a morphism of inverse systems in the category \
$\Cal {CGR}$.

Choose a point \ $\widetilde{e} \in p^{-1}(e)$  \.
Let $\widetilde e_{\lambda}:=h_{\lambda}({\widetilde e})$ \
for each index $\lambda\in\Lambda_\beta$.
It follows from commutativity of the diagram (4) that \
$\widetilde e_\lambda\in p^{-1}_{\lambda}(e_\lambda)$, where \
$e_\lambda$ \ is the identity of the group \ $G_\lambda$.

Applying Theorem to the covering mappings \
$p_{\lambda}: X_{\lambda}\to G_{\lambda}, \lambda \in \Lambda_\beta,$
we endow each space \ $X_\lambda$ \ with the structure of topological group
such that \ $\widetilde e_\lambda$ \ is the identity and the $k$-fold covering
mapping \ $p_{\lambda}$ \ is a homomorphism of groups.

We need to show that all mappings $h_{\lambda}^{\mu}: X_\mu\to X_\lambda$
become homomorphisms of groups. For this purpose, we fix a mapping
$h_{\lambda}^{\mu}$ and consider the diagram

$$
\CD
                 @.      X_\lambda  \\
@.                     @VV{p_\lambda}V  \\
X_\mu\times X_\mu    @>F>>  G_\lambda
\endCD
$$
where the mapping \ $F$ \ is defined by setting \
$F(x,y)= p_\lambda{\circ}h_{\lambda}^{\mu}(xy)$.
Denote by \ $F_1$ \ and \ $F_2$ \ the mappings from the Cartesian product \
$X_\mu\times X_\mu$ \ to the space \ $X_\lambda$ \ given by
$$F_1(x,y)=h_{\lambda}^{\mu}(xy) \quad \text{and} \quad
F_2(x,y)=h_{\lambda}^{\mu}(x)h_{\lambda}^{\mu}(y), \quad x,y\in X_\mu.$$
Both these mappings are liftings of the mapping \ $F$ \ to \ $X_\lambda$.
That is,
making use of \ $F_1$ \ or \ $F_2$, the above diagram can be completed
to a commutative one: $p_\lambda \circ F_1=F $ and $p_\lambda \circ F_2=F $.
Indeed, since the diagram (10) is commutative and the
mappings $p_\lambda$, $p_\mu$ and $\pi_{\lambda}^{\mu}$
are homomorphisms of groups, we have:
$${p_\lambda\circ}h_{\lambda}^{\mu}(xy) =
p_{\lambda}(h_{\lambda}^{\mu}(x)h_{\lambda}^{\mu}(y))
\quad \text{for any}\quad x,y\in X_\mu.$$
In addition, $F_1(\widetilde e_\mu,\widetilde e_\mu)=
F_2(\widetilde e_\mu,\widetilde e_\mu),$ where
$\widetilde e_\mu$ is the identity of the group $X_\mu$.
Because the space $X_\mu\times X_\mu$  is connected,
it follows from the uniqueness property for liftings to covering spaces
\cite{18; Chapter 2, \S 2, Theorem 2} that
$F_1=F_2$. This means that \ $h_{\lambda}^{\mu}: X_\mu\to X_\lambda$ \ is
a homomorphism of groups.

Thus the family
$$\{\,p_\lambda: \lambda\in\Lambda_\beta\,\}:
\{\,X_\lambda,h_{\lambda}^{\mu},\Lambda_\beta\,\} \to
\{\,G_\lambda,\pi_{\lambda}^{\mu},\Lambda_\beta\,\}$$
becomes a morphism of inverse systems in the category $\Cal {CGR}$.
Therefore \linebreak $p_\infty:X_\infty\to G_\infty$ \ is a homomorphism of
compact groups. If the group $G$ is abelian, then each $k$-fold covering
mapping \ $p_\lambda:X_\lambda\to G_\lambda$ \ is a homomorphism between
abelian groups, so that \ $p_\infty:X_\infty\to G_\infty$ \ is a homomorphism
of abelian groups as well.

Finally, it remains to note that now in the diagram (11) the mappings
$\sigma$ and $p_\infty$ are homomorphisms of compact groups.
Using this diagram, one can easily equip the space $X$ with the desired
structure of the topological group whose identity is the point
$\widetilde{e}$.
The uniqueness of such structure follows from the uniqueness property
for liftings to covering spaces.\qed
\enddemo

Note that Theorem~1 was proved in \cite{19} for
a finite-sheeted covering mapping onto a solenoidal group by a different
method.

The following examples show that, in general, the covering group theorem fails
for a finite-sheeted covering mapping with a non-connected covering space.
As usual, in all examples the symbol $\Bbb S^1$ stands for the unit circle.

\example{Example 1} Let $X_k=\Bbb S^1\times\{k\}$, where $k=1,2$.
Denote by $p_k:X_k\to \Bbb S^1$ the~$k$-fold
covering mapping sending $(z,k)$ to $z^k$, $z\in \Bbb S^1$.
Consider the combination $p_1{\bigtriangledown}p_2$ of the mappings
$p_1$ and $p_2$ which is defined as follows:
$$p_1{\bigtriangledown}p_2: X_1{\oplus}X_2\to \Bbb S^1:x\mapsto p_k(x)
\quad \text{if} \quad x\in X_k.$$
Here, and in what follows, the symbol ${\oplus}$ stands for the sum of
topological spaces.
We denote the space \ $X_1{\oplus}X_2$ \ and \ the mapping \
$p_1{\bigtriangledown}p_2$ \
by \ $Y$ \ and \ $p$ \ respectively. Clearly, the mapping \ $p:Y\to \Bbb S^1$ \
is the three-fold covering mapping.

We claim that the structure of topological group on \ $\Bbb S^1$ \ does not lift
to the~covering space \ $Y$.
To obtain a contradiction, we suppose that there exists a structure
of topological group on \ $Y$ \ turning \ $p$ \ into a homomorphism of groups.
Let us assume that the identity of the group \ $Y$ \ belongs to \ $X_2$.
One can easily see that for each \ $y\in X_1$ \ the equality
\ $X_1=yX_2$ \ holds.  This implies that the~restrictions
\ $p\vert_{X_1}=p_1$ \ and \ $p\vert_{X_2}=p_2$ \ have the same degree.
The case when the~identity of \ $Y$ \ lies in \ $X_1$ \ is similar.
This contradiction proves our claim.
\endexample

Using the covering mapping \ $p:Y\to \Bbb S^1$ \ from Example 1,  we shall
construct a finite-sheeted covering mapping onto a non-connected compact
group.

\example{Example 2} For $j=0,1$ we let  $Y_j=Y \times \{j\}$.
Denote by $q_0$ and $q_1$ the mappings defined by the formulae:
$$q_j:Y_j \to \Bbb S^1 \times \{j\}:(y,j) \mapsto (p(y),j),
\quad y \in Y, \quad j=0,1.$$
The sum $q_0{\oplus}q_1: Y_0{\bigoplus}Y_1\to \Bbb S^1 \times \Bbb Z_2$ of
the mappings $q_0$ and $q_1$ which is given by
$$
q_0{\oplus}q_1(y)=q_j(y)\quad \text{if} \quad y \in Y_j, \quad j=0,1,
$$
is the required covering mapping.
Here \ $\Bbb Z_2=\{0,1\}$ \ is the group of order two equipped with the discrete
topology, and \ $\Bbb S^1 \times \Bbb Z_2$ \ is the direct product of groups.
\endexample

\remark{Remark}
We refer the reader to \cite{20} and \cite{21} for conditions for
the lifting of topological group structure on a non-connected topological
group to a covering space.
\endremark

Theorem 1 suggests the following question.
\proclaim{Question 1} Can one take the group $G$ to be locally compact
in Theorem 1?
\endproclaim

In connection with the above theorems the following interesting
conjecture was stated by Professor S. A. Bogatyi.

\proclaim{Conjecture} The covering group theorem holds for an overlay
mapping from a~connected space onto a topological group.
\endproclaim

\head
4. Covering mappings onto abelian groups
\endhead

In this section \ $G$ \ is a compact connected abelian group, and \
$\widehat{G}$ \ is its additive character group. We refer to the appropriate
sections of \cite{4} and \cite{22} for the fundamental facts about
the character group and Pontryagin duality. If \ $\pi:H_1\to H_2$ \
is a homomorphism between two compact abelian topological groups,
then \ $\widehat{\pi}:\widehat{H_2}\to \widehat{H_1}$ \ denotes the
dual homomorphism between the character groups, i.e, \
$\widehat{\pi}(\chi)=\chi \circ \pi$ \ for \ $\chi \in \widehat{H_2}$.

For given \ $k\in \Bbb N$ \ one has  the homomorphism
$$\tau_k :\widehat G\to\widehat G:\chi \mapsto k\chi.$$
If \ $\tau_k$ \ is surjective, then the group \ $\widehat G$ \
is said to be $k$-{\it divisible}, or we say that the~group \ $\widehat G$ \
{\it admits division by} $k$. Note that if \ $\widehat G$ \ is $k$-divisible,
then the~homomorphism \ $\tau_k$ \ is an automorphism since
\ $\widehat G$ \ is a torsion free abelian group (\,\cite{22; Theorem 24.25}\,).
Necessary and sufficient conditions for
\ $\widehat G$ \ to be $k$-divisible are given in \cite{23}. In particular, it is shown
\cite{23; Theorem 1.1}
that the character group \ $\widehat G$ \ is $k$-divisible if and only if
the group $G$ \ admits a $k$-mean. 
There is a large literature
concerning the problem on the existence of means on topological spaces (see,
e.g., \cite{23} -- \cite{28} where additional references are contained).
Recall that, for $k\in \Bbb N$, a~$k$-{\it mean} on \ $G$ \ is
a mapping \ $\mu:G\times G \times \ldots \times G \to G$ \ from
the Cartesian product of $k$ copies of \ $G$ \ such that \
$\mu(g_1,g_2,\ldots,g_k)=
\mu(g_{\pi(1)},g_{\pi(2)},\ldots,g_{\pi(k)})$ \ for any permutation \ $\pi$ \
of the set \ $\overline{k}$ \ and \ $\mu(g,g,\ldots,g)=g$ \ for all $g\in G$.

As a corollary of Theorem~1, we obtain the following theorem concerning
finite-sheeted covering mappings onto abelian groups.

\proclaim{Theorem 2}
Let \ $p:X \to G$ \ be a finite-sheeted covering mapping 
from a connected topological space \ $X$ \ onto a compact connected abelian
group \ $G$ \ such that the~character group \ $\widehat G$ \ admits division
by the degree of \ $p$.
Then \ $p$ \ is a~homeomorphism.
\endproclaim

\demo{Proof}
According to Theorem~1, we can consider the mapping \ $p$ \
as a homomorphism between compact connected abelian groups.
Let \ $e$ \ and \ $\widetilde{e}$ be identities of the groups \ $G$ \ and
\ $X$ \ respectively. Denote by \ $\tau: X \to X/p^{-1}(e)$ \
the canonical homomorphism onto the factor group \ $X/p^{-1}(e)$ \ and
by \ $\sigma: X/p^{-1}(e) \to G$ \ the~isomorphism of topological groups
such that \ $p= \sigma \circ \tau$.

We claim that \ $\tau$ \ is an isomorphism of topological groups.
In order to prove this assertion we shall show that \
$p^{-1}(e)=\{\widetilde{e}\}$.

To obtain a contradiction we suppose that \ $p^{-1}(e) \neq \{\widetilde{e}\}$.
Obviously, in this case, there exists a character \ $\chi\in\widehat X$ \
such that the restriction of \ $\chi$ \ to the subgroup \ $p^{-1}(e)$ \
is not the identity character.

Denote by \ $k$ \ the degree of \ $p$.
Since \ $p^{-1}(e)$ \ is a finite group of order \ $k$, the restriction
of the character \ $k\chi$ \ to \ $p^{-1}(e)$ \ is the identity character.
As a~consequence \cite{22, Theorem 23.25}, there is a character \
$\theta\in\widehat{X/p^{-1}(e)}$ \ such that
$$\widehat\tau(\theta) = k\chi.$$

The group \ $\widehat{X/p^{-1}(e)}$ \ is \ $k$-divisible. Hence there
is a character \ $\eta\in\widehat{X/p^{-1}(e)}$  \ satisfying the equality
$$k\eta =\theta,$$
which implies that $$k\widehat{\tau}(\eta) = k\chi.$$
But \ $\widehat X$ \ is a torsion free group \cite{22, Theorem 24.25}, so that
$$\widehat{\tau}(\eta) = \chi.$$

It follows that the restriction of $\chi$ to the
subgroup \ $p^{-1}(e)$ \ is the identity character. This contradiction yields
the desired equality $p^{-1}(e)=\{\widetilde{e}\}$.

Thus \ $\tau$ \ and \ $p$ \ are isomorphisms of topological groups,
as required.\qed
\enddemo

As immediate consequences of Theorem~2 and Theorem~1.1 from \cite{23}
we have the following statements.

\proclaim{Corollary 1}
Let \ $k \geq 2$.  If the character group \ $\widehat G$ \ is
$k$-divisible, then there exists no a $k$-fold covering mapping from a connected
topological space onto \ $G$.
\endproclaim

\proclaim{Corollary 2}
Let \ $k \geq 2$. If there is a $k$-fold covering mapping from a connected
topological space onto \ $G$, then \ $G$ \ does not admit a $k$-mean.
\endproclaim

\example{Example 3}
Let $P=(\, p_1, p_2, \ldots)$ be a sequence of prime numbers (1 not being
included as a prime) and let \ $\Sigma_P$ \ be the {\it $P$-adic solenoid},
that is, the inverse limit of the inverse sequence
(= the inverse system over $\Bbb N$)
$$
\CD
\Bbb S^1 @< \pi_1^2 << \Bbb S^1 @< \pi_2^3 << \Bbb S^1 @< \pi_3^4 << \cdots ,
\endCD
$$
where the bonding mappings are given by \
$\pi_n^{n+1}(z)=z^{p_n}$ for all $z \in \Bbb S^1$ and $n \in \Bbb N$.
As is well known, the $P$-adic solenoid is a compact connected abelian group.

The character group of \ $\Sigma_P$ \ is isomorphic to the discrete
additive group
of all rationals of the form \ $\frac{m}{p_1p_2\ldots p_n}$,
where \ $n\in \Bbb N$ \ and \ $m$ is an integer \cite{22; (25.3)}.

Assume that \ $k$ \ is a positive integer such that each its prime divisor
is equal to infinitely many terms of the sequence $P$.
By number-theoretic considerations one can easily see that
the group \ $\widehat{\Sigma_P}$ \ is $k$-divisible. According to Theorem~1,
the~$P$-adic solenoid \ $\Sigma_P$ \ does not admit \ $k$-fold covering
mapping from a connected topological space onto itself.

If for each prime number \ $p\in \Bbb N$ \
we have \ $p=p_n$ \ for infinitely many indices \ $n$, then
the group \ $\widehat{\Sigma_{P}}$ \ is $k$-divisible for any $k\in \Bbb N$.
In this case the $P$-adic solenoid \ $\Sigma_{P}$ \ admits only trivial
finite-sheeted covering mappings from topological spaces onto itself.
\endexample

\remark{Remark}
One can construct (\, see \cite{6, Example 2}, \cite{8, Proposition~2.2}\,)
an example of solenoid \ $\Sigma_{P}$ \ which
does not admit a $k$-mean and such that there is no a $k$-fold covering
mapping from a connected space onto \ $\Sigma_{P}$.
\endremark

Using Theorem~2, we obtain the following criterion of triviality of
finite-sheeted covering mappings onto compact connected abelian groups.

\proclaim{Theorem 3}
All finite-sheeted covering mappings of degree $k\in \Bbb N$ onto a compact
connected abelian group $G$ are trivial if and only if the character group \
$\widehat G$ \ is $k!$-divisible.
\endproclaim

\demo{Proof} As was mentioned in Introduction, we need only to prove
the sufficiency.

Suppose that the group \ $\widehat G$ \ is $k!$-divisible.
Let \ $p: Y\to G$ \ be a $k$-fold covering mapping, where \ $Y$ \
is a non-connected topological space.
Choose a component \ $Y_0$ \ of the space \ $Y$ \  and
consider a finite-sheeted covering mapping $p\vert_{Y_0} : Y_0\to G$, i.e.,
the restriction of $p$ to $Y_0$. Obviously, the group \ $\widehat G$ \
admits division by degree of the covering mapping \ $p\vert_{Y_0}$.
Applying Theorem~2 to \ $p\vert_{Y_0}$, we conclude that
the~mapping \ $p\vert_{Y_0}$ is a homeomorphism. Now one can easily show
that \ $p$ \ is a~trivial $k$-fold covering mapping.\qed
\enddemo

\remark{Remark} As an application of Theorem~1, it can be shown \cite{29} that each
finite-sheeted covering mapping from (in general, non-connected) topological
space onto \ $G$ \ is equivalent to a covering mapping which is generated
by a separable polynomial over the Banach algebra \ $C(G)$. This fact
together with the result by E.~A.~Gorin and V.~Ya.~Lin mentioned in
Introduction yields another proof of Theorem~3.
\endremark

In connection with Theorem~3, an interesting question is the following.

\proclaim{Question 2} Is it true that all $k$-fold covering mappings onto
a compact connected topological space \ $X$ \ are trivial provided that
the first \v{C}ech cohomology group of \ $X$ \ with integer
coefficients \ $H^1(X, \Bbb Z)$ \ is \ $k!$-divisible?
\endproclaim

Finally, combining Theorem 3 with Theorem 1.1 from \cite{23},  we have:

\proclaim{Theorem 4}
All finite-sheeted covering mappings of degree $k\in \Bbb N$ onto a compact
connected abelian group $G$ are trivial if and only if \ $G$ \
admits a $k!$-mean.
\endproclaim


\Refs\nofrills{References}
\ref \no 1
\by Gorin E. A. and Lin V. Ya.
\paper Algebraic equations with continuous coefficients and some problems
of the algebraic theory of braids
\jour Matem. Sbornik. 78 (120) (1969)
\pages  579--610  (Russian)
\endref
\ref \no 2
\by Gorin E. A. and Lin V. Ya.
\paper On separable polynomials over commutative Banach algebras \linebreak
\jour Doklady Akad. Nauk SSSR. 218 (1974)
\pages  505--508 (Russian)
\endref
\ref \no 3
\by Zyuzin Yu. V.
\paper Irreducible separable polynomials with holomorphic coefficients
on a certain class of complex spaces
\jour  Matem. Sbornik. 102 (144) (1977)
\pages  569--591 (Russian)
\endref
\ref \no 4
\by Pontryagin L. S.
\book Continuous Groups
\publ Nauka, Moscow
\yr 1984  \ (Russian)
\endref
\ref \no 5
\by  Massey W. S.
\book Algebraic topology: An Introduction
\publ Harcourt Brace Jovanovich, New York
\yr 1967
\endref
\ref \no 6
\by Fox R.H.
\paper On shape
\jour Fund. Math. 74 (1972)
\pages  47 -- 71
\endref
\ref \no 7
\by Fox R.H.
\book Shape theory and covering spaces
\publ Lecture Notes in Math. 375, 71--90, Topology Conference Virginia
Polytechnic Institute, 1973 (R.F. Dickman Jr., P.Fletcher, eds.), Springer,
Berlin --- Heidelberg --- New York
\yr 1974
\endref
\ref \no 8
\by Moore T.T.
\paper On Fox's theory of overlays
\jour Fund. Math. 99 (1978)
\pages  205 -- 211
\endref
\ref \no 9
\by Ko\c{c}ak \c{S}.
\paper On the fundamental theorem of overlays
\jour Note Mat. 10 (1990)
\pages  355 -- 362
\endref
\ref \no 10
\by Mrozik P.
\paper Images and exactness in the category of pro-groups and the lifting
problem for covering projections
\jour Glasnik Mat. 28 (1993)
\pages  209 -- 226
\endref
\ref \no 11
\by Marde\v{s}i\'{c} S. and Matijevi\'{c} V. 
\paper Classifying overlay structures of topological spaces
\jour Topology Appl. 113 (2001)
\pages 167 -- 209
\endref
\ref \no 12
\by Grigorian S. A. and Gumerov R. N.
\book A criterion of triviality for finite-sheeted coverings of compact
connected abelian groups
\publ Trudy Lobachevskii Matem. Centre, V. 5, 238--240,
Proceedings of International Conference, Unipress, Kazan
\yr 2000
\endref

\ref \no 13
\by Grigorian S.A. and Gumerov R.N.
\paper On a covering group theorem and its applications
\jour Loba\-chevskii J. Math. X (2002)
\pages  9--16
\endref

\ref \no 14
\by  Eilenberg S. and Steenrod N.
\book Foundations of Algebraic topology
\publ Princeton Univ. Press, Princeton, N. J.
\yr 1952
\endref

\ref \no 15
\by Engelking R.
\book General topology
\publ Monografie Matematyczne, Vol. 60, Polish Scientific Publishers, Warszawa
\yr 1977
\endref

\ref \no 16
\by Bucur I. and Deleanu A., with the collaboration of Hilton P.J.
\book Introduction to the Theory of Categories and Functors.
\publ Pure and Applied Mathematics, Vol. XIX, Wiley - Interscience Publ.,
London - New York - Sydney
\yr 1968
\endref
\ref \no 17
\by Weil A.
\book L'Integration dans les Groupes Topologiques et ses Applications.
\publ Publ. L'Institut Math. de Clermont-Ferrand, Paris
\yr 1940
\endref
\ref \no 18
\by Spanier E. H. \nofrills
\book Algebraic Topology. \nofrills
\publ McGraw-Hill, New York
\yr 1966
\endref
\ref \no 19
\by Grigorian S.A., Gumerov R.N. and Kazantsev A.V. 
\paper Group structure in finite coverings of compact solenoidal groups
\jour Lobachevskii J. Math. VI (2000)
\pages 39--46
\endref

\ref \no 20
\by Taylor R. L.
\paper Covering groups of non-connected topological groups
\jour Proc. Amer. Math. Soc. 5 (1954)
\pages 753--768
\endref
\ref \no 21
\by Brown R. and Mucuk O.
\paper Covering groups of non-connected topological groups revisited
\jour Math. Proc. Cambridge Phil. Soc. 115 (1994)
\pages 97--110
\endref

\ref \no 22
\by Hewitt E. and Ross K. A.
\book Abstract harmonic analysis, Vol. I
\publ Springer-Verlag, Berlin
\yr 1963
\endref

\ref \no 23
\by Keesling J.
\paper The group of homeomorphisms of a solenoid
\jour Trans. Amer. Math. Soc. 172 (1972)
\pages 119--131
\endref

\ref \no 24
\by Eckmann B.
\paper R\"{a}ume mit Mittelbildungen
\jour Comment. Math. Helv. 28 (1954)
\pages 329--340
\endref

\ref \no 25
\by Bacon P.
\paper Unicoherence in means
\jour Colloq. Math. 21 (1970)
\pages 211--215
\endref


\ref \no 26
\by Charatonik J. J.
\book Some problems concerning means on topological spaces
\publ Topology, Measures and Fractals (C.Bandt, J. Flachsmayer and H. Haase,
eds.), Mathematical Research, vol. 66, 166--177, Akademie Verlag, Berlin
\yr 1992
\endref

\ref \no 27
\by Kawamura K. and Tymchatyn E. D.
\paper Continua which admit no mean
\jour Colloq. Math. 71 (1996)
\pages 97--105
\endref

\ref \no 28
\by  Illanes A. and Nadler Jr., S. B.
\book Hyperspaces
\publ Marcel Dekker, Inc., New York
\yr 1999
\endref

\ref \no 29
\by Grigorian S.A. and Gumerov R.N.
\paper On polynomials with continuous coefficients and finite coverings
of compact connected abelian groups \nofrills
\jour (manuscript)
\endref


\endRefs
\enddocument